\documentclass[12pt]{article}                                               

\input{amssym.tex}                                                           
                                                                             
\begin{document}                                                             
\title{On a correlation between ranks of elliptic curves and periods of continued fractions}

\author{Igor  ~Nikolaev}


\date{}
 \maketitle


\newtheorem{thm}{Theorem}
\newtheorem{lem}{Lemma}
\newtheorem{dfn}{Definition}
\newtheorem{rmk}{Remark}
\newtheorem{cor}{Corollary}
\newtheorem{prp}{Proposition}
\newtheorem{exm}{Example}
\newtheorem{cnj}{Conjecture}

\newcommand{\ch}{\hbox{\bf Char}}
\newcommand{\jac}{\hbox{\bf Jac}}
\newcommand{\ka}{\hbox{\bf k}}
\newcommand{\n}{\hbox{\bf n}}
\newcommand{\mod}{\hbox{\bf mod}}
\newcommand{\Z}{\hbox{\bf Z}}
\newcommand{\Q}{\hbox{\bf Q}}

\newcommand{\Coh}{\hbox{\bf Coh}}
\newcommand{\Mod}{\hbox{\bf Mod}}
\newcommand{\Irred}{\hbox{\bf Irred}}
\newcommand{\Spec}{\hbox{\bf Spec}}
\newcommand{\Tors}{\hbox{\bf Tors}}

\begin{abstract}
It is proved that the rank of elliptic curves  with complex multiplication 
 introduced  by B.~H.~Gross 
is  one less an arithmetic complexity of the  corresponding 
noncommutative tori with real multiplication.

\vspace{7mm}


\vspace{5mm}
\end{abstract}

\section{Introduction}
Let ${\cal E}(k)$ be a non-singular  elliptic curve defined  over a field 
$k\subseteq {\Bbb C}$. 
Recall that the  Sklyanin algebra $S(\alpha,\beta,\gamma)$ is  a  free $k$-algebra on four generators $\{x_1,\dots, x_4\}$   
and six  quadratic relations:  
\begin{equation}
\left\{
\begin{array}{ccc}
x_1x_2-x_2x_1 &=& \alpha(x_3x_4+x_4x_3),\\
x_1x_2+x_2x_1 &=& x_3x_4-x_4x_3,\\
x_1x_3-x_3x_1 &=& \beta(x_4x_2+x_2x_4),\\
x_1x_3+x_3x_1 &=& x_4x_2-x_2x_4,\\
x_1x_4-x_4x_1 &=& \gamma(x_2x_3+x_3x_2),\\ 
x_1x_4+x_4x_1 &=& x_2x_3-x_3x_2,
\end{array}
\right.
\end{equation}
where $\alpha,\beta,\gamma\in k$ and  
$\alpha+\beta+\gamma+\alpha\beta\gamma=0$  \cite[Example 8.5]{StaVdb1}.  The algebra $S(\alpha,\beta,\gamma)$ is 
a twisted homogeneous coordinate ring  of the elliptic curve ${\cal E}(k)\subset {\Bbb C}P^3$ 
 given in the  Jacobi form $u^2+v^2+w^2+z^2 =
{1-\alpha\over 1+\beta}v^2+{1+\alpha\over 1-\gamma}w^2+z^2 = 0$;
we refer the reader to \cite{StaVdb1} for the missing definitions and details.

 Consider a self-adjoint representation 
$\rho:  S(\alpha,\beta,\gamma)\to {\cal B}({\cal H})$,  where ${\cal B}({\cal H})$ is the ring of 
bounded linear operators on a Hilbert   space ${\cal H}$.   The norm-closure of  $\rho(S(\alpha,\beta,\gamma))$ 
is a $C^*$-algebra \cite{M}. The $C^*$-algebra is isomorphic to the so-called noncommutative torus
${\cal A}_{\theta}$,    i.e.   a $C^*$-algebra  generated by two unitary operators $u$ and $v$
satisfying the commutation relation $vu=e^{2\pi i\theta}uv$  for a real constant  $\theta$ \cite[Section 1.3]{N} and \cite{Rie2}.
The map ${\cal E}(k)\mapsto {\cal A}_{\theta}$ is a functor, 
such that if  the curves ${\cal E}(k)$ and ${\cal E}'(k)$  are isomorphic over $k$,  then  the 
algebras ${\cal A}_{\theta}$ and ${\cal A}_{\theta'}$ are isomorphic and   if  ${\cal E}(k)$ and ${\cal E}'(k)$ are isomorphic over ${\Bbb C}$,  then 
${\cal A}_{\theta}$ and ${\cal A}_{\theta'}$ are Morita equivalent,  i.e.  ${\cal A}_{\theta}\otimes {\cal K}\cong
{\cal A}_{\theta'}\otimes {\cal K}$,  where ${\cal K}$ is the $C^*$-algebra of compact operators \cite[Corollary 1.2]{Nik3}.  
If $\theta$ is a quadratic irrationality,  the  algebra $ {\cal A}_{\theta}$
is said to have  real multiplication \cite{Man1}.

Denote by   ${\cal E}_{CM}^{(-D,f)}$
an  elliptic curve with  complex multiplication by an order $R_f$ of conductor $f$
in the imaginary quadratic field  ${\Bbb Q}(\sqrt{-D})$
\cite[pp. 95-96]{S}.  
The noncommutative torus corresponding to ${\cal E}_{CM}^{(-D,f)}$
has  real multiplication by  an order ${\goth R}_{\goth f}$ of 
conductor ${\goth f}$  in the quadratic field ${\Bbb Q}(\sqrt{D})$;  such a torus
we denote by  ${\cal A}_{RM}^{(D, {\goth f})}$.    The conductor ${\goth f}$ 
is defined from the equation $|Cl~({\goth R}_{\goth f})|=|Cl~(R_f)|$,  where $Cl$ is the
class  group of the respective orders \cite[Theorem 6.1.3]{N}.

Let $({\cal E}_{CM}^{(-D,f)})^{\sigma},  ~\sigma\in Gal~(k|{\Bbb Q})$ be 
the  Galois conjugate of the curve  ${\cal E}_{CM}^{(-D,f)}$; 
by a {\it ${\Bbb Q}$-curve} one understands   ${\cal E}_{CM}^{(-D,f)}$,
such that  there exists  an  isogeny between $({\cal E}_{CM}^{(-D,f)})^{\sigma}$ and   ${\cal E}_{CM}^{(-D,f)}$
for  each    $\sigma\in Gal~(k|{\Bbb Q})$. 
Let ${\goth P}_{3 ~\mod ~4}$ be the set of all primes  
$p=3 ~\mod ~4$;    it is known that  ${\cal E}_{CM}^{(-p,1)}$ is a ${\Bbb Q}$-curve 
whenever $p\in  {\goth P}_{3 ~\mod ~4}$ \cite[p. 33]{G}. 
The rank of  ${\cal E}_{CM}^{(-p,1)}$  is always divisible by $2h_K$,  where $h_K$ is the 
 class number of field $K:={\Bbb Q}(\sqrt{-p})$  \cite[p.49]{G};  
  by a {\it ${\Bbb Q}$-rank}  of  ${\cal E}_{CM}^{(-p,1)}$
 we understand the integer 
 $rk_{\Bbb Q}({\cal E}_{CM}^{(-p,1)}):={1\over 2h_K}~rk~({\cal E}_{CM}^{(-p,1)})$.

Denote by $(\overline{a_1,\dots, a_P})$ the minimal period of 
continued fraction of $\sqrt{D}$;  by an  {\it arithmetic complexity}  of 
the algebra  ${\cal A}_{RM}^{(D, {\goth f})}$
we understand  the   number of independent   $a_i$  in  
the period $(\overline{a_1,\dots, a_P})$,  see Section 2
for an exact definition.
 The complexity is   denoted by  $c({\cal A}_{RM}^{(D, {\goth f})})$ and it
 is equal to the Krull dimension of  connected component of 
an affine variety given by the diophantine equation (\ref{eq5}).

Since the algebra ${\cal A}_{RM}^{(D, {\goth f})}$  encodes geometry of the curve ${\cal E}_{CM}^{(-D,f)}$,
it is natural to expect that the values of $c({\cal A}_{RM}^{(D, {\goth f})})$ and $rk_{\Bbb Q}({\cal E}_{CM}^{(-D,f)})$
are related.  (For simplicity, we further assume that $D=p$ is a prime number and  $f={\goth f}=1$.)     The aim of our note is a formula
linking the ${\Bbb Q}$-ranks of the ${\cal E}_{CM}^{(-p,1)}$ with the arithmetic complexity of the 
algebra ${\cal A}_{RM}^{(p, 1)}$. 
\begin{thm}\label{thm1}
$rk_{\Bbb Q}~({\cal E}_{CM}^{(-p,1)}) +1=c({\cal A}_{RM}^{(p,1)})$,
whenever $p\equiv 3 ~\mod ~4$. 
\end{thm}
\begin{rmk}
\textnormal{
It is known that the integer $rk_{\Bbb Q}~({\cal E}_{CM}^{(-p,1)})$ depends on a twist 
of the ${\cal E}_{CM}^{(-p,1)}$. On the other hand,  the number $c({\cal A}_{RM}^{(p,1)})$ depends
on the period
$(\overline{a_1,\dots, a_P})$ alone. This observation does not contradict
theorem \ref{thm1},  since it is known that the twists of  ${\cal E}_{CM}^{(-p,1)}$ correspond to
continued fractions  of the form
$[a_0,\dots, a_k; \overline{a_{k+1},\dots, a_{k+P}}]$. 
The arithmetic complexity of  such fractions is defined in \cite[Section 6.2.1]{N}.
In particular,  the ${\cal E}_{CM}^{(-p,1)}$ corresponds to a  continued fraction
$[a_0; \overline{a_{1},\dots, a_{P}}]$,  see Perron's Lemma \ref{lm1}. 
 Thus  the formula $rk_{\Bbb Q}~({\cal E}_{CM}^{(-p,1)}) +1=c({\cal A}_{RM}^{(p,1)})$
is correctly defined. 
}  
\end{rmk}
\begin{rmk}
\textnormal{
A generalization of Theorem \ref{thm1} is proved in  \cite[Theorem 6.2.1]{N} by different methods.
The value of present paper is a direct approach using  the Perron  Lemma \cite[p. 88]{P} 
and  construction of  explicit  examples based on the Gross  Thesis  \cite[p. 78]{G}. 
}  
\end{rmk}
The article is organized as follows.  The arithmetic complexity is 
defined  in Section 2. Theorem \ref{thm1}  is proved in Section 3.
In Section 4 we illustrate theorem \ref{thm1} by examples of ${\cal E}_{CM}^{(-p,1)}$ 
 for  primes $p<100$.

\section{Arithmetic complexity}
Let $[a_0,a_1,\dots]$ be a regular continued fraction and consider 
  a system of the linear equations:
\begin{equation}\label{eq2}
\left\{
\begin{array}{ccc}
y_0 &=& a_0y_1+y_2\\
y_1 &=& a_1y_2+y_3\\
y_2 &=& a_2y_3+y_4\\
       &\vdots& 
\end{array}
\right.
\end{equation}
One can  write  (\ref{eq2}) in the form:
\begin{equation}\label{eq3}
\left\{
\begin{array}{ccc}
y_j &=& A_{i-1,j}y_{i+j}+a_{i+j}A_{i-2,j}y_{i+j+1}\\
y_{j+1} &=& B_{i-1,j}y_{i+j}+a_{i+j}B_{i-2,j}y_{i+j+1},
\end{array}
\right.
\end{equation}
where the polynomials $A_{i,j}$ and $B_{i,j}$   are called continuants 
(or, Muir's symbols) \cite[p. 10]{P}.  
\begin{lem}\label{lm1}
{\bf  (\cite[p. 88  and p. 107]{P})} 
There exists  a square-free integer $D>0$,  such that
\displaymath
 [x_0; \overline{x_1, x_2, \dots, x_2, x_1, x_P}]=
 \cases{\sqrt{D}, & if  $x_P=2x_0$  and $D=2,3~\mod ~4,$\cr
                {\sqrt{D}+1\over 2}, & if   $x_P=2x_0-1$ and  $D=1 ~\mod ~4,$}
\enddisplaymath
if and only if   $x_P$ satisfies  an Euler equation:
\begin{equation}\label{eq5}
x_P=mA_{P-2, 1}-(-1)^PA_{P-3,1}B_{P-3,1},
\end{equation}
where  $m>0$ is an integer.  
\end{lem}
Denote by {\bf x}$_D=(x_0,\dots, x_P)$ a solution of the Euler equation (\ref{eq5})
and consider an affine algebraic set ${\cal A}$ defined by the polynomial 
equation (\ref{eq5}). 
By  an  {\it  Euler variety}   $V_D$ one  understands  the projective closure
of an irreducible component of  ${\cal A}$ containing the point  {\bf x}$_D$. 
\begin{dfn}
By an arithmetic complexity $c({\cal A}_{RM}^{(D,1)})$ of the algebra
 ${\cal A}_{RM}^{(D,1)}$  we  understand  the Krull  dimension of 
 the Euler variety $V_D$. 
  \end{dfn}
\begin{rmk}
\textnormal{
Roughly speaking, the arithmetic complexity is the number of independent
$x_i$ in  the continued fraction $[x_0; \overline{x_1, x_2, \dots, x_2, x_1, x_P}]$.
In particular,   for  the ``generic''  quadratic irrationalities $\{\theta=r_1+r_2\sqrt{D} ~|~r_1,  r_2\in {\Bbb Q}\}$
the arithmetic complexity is equal to $P$. 
}  
\end{rmk}
\begin{exm}\label{ex1}
{\bf (\cite[p. 90]{P})}
{\normalfont
If  $P=4$,   then the continuants  are: $A_{P-3,1}=A_{1,1}=x_1x_2+1$, ~$B_{P-3,1}=B_{1,1}=x_2$ and 
$A_{P-2,1}=A_{2,1}=x_1x_2x_3+x_1+x_3=x_1^2x_2+2x_1$, since $x_3=x_1$.  Thus, 
the Euler equation (\ref{eq5}) takes the form: 
\begin{equation}\label{eq6}
2x_0= m(x_1^2x_2+2x_1)-x_2(x_1x_2+1),
\end{equation}
and, therefore,  $\sqrt{x_0^2+m(x_1x_2+1)-x_2^2}=[x_0,\overline{x_1,x_2,x_1,2x_0}]$. 
Let us show that the affine set defined by (\ref{eq6}) is reducible. 
Indeed,  by lemma \ref{lm1},  parameter $m$ must be integer for all  integer values
of $x_0,x_1$ and $x_2$.  This is not possible in general, since from (\ref{eq6})
one obtains  $m=(2x_0+x_2(x_1x_2+1))(x_1^2x_2+2x_1)^{-1}$ is a rational number. 
However,  a restriction to  $x_1=1, ~x_2=x_0-1$ defines a
connected component of the affine set (\ref{eq6}),  since in this case $m=x_0$ is always an integer.    
Thus,  one gets a family of solutions of (\ref{eq6}) of the form
 $\sqrt{(x_0+1)^2-2}=[x_0,\overline{1, x_0-1, 1, 2x_0}]$. 
 We conclude that $c({\cal A}_{RM}^{(D,1)})=1$,
 where $D=(x_0+1)^2-2$.
  }
 \end{exm}

\section{Proof of theorem  \ref{thm1}}
We shall split the proof in a series of lemmas. 
\begin{lem}\label{lm3}
If $[x_0,\overline{x_1,\dots, x_k,\dots, x_1, 2x_0}]\in\sqrt{{\goth P}_{3~\mod~4}}$  ,
then:

\smallskip
(i) $P=2k$ is an even number,  such that:

\smallskip
\hskip1cm
(a) $P\equiv 2~\mod~4$,  if $p\equiv 3~\mod~8$;

\smallskip
\hskip1cm
(b) $P\equiv 0~\mod~4$,  if $p\equiv 7~\mod~8$;

\medskip
(ii) either of two is true:

\smallskip
\hskip1cm
(a) $x_k=x_0$ (a culminating period);

\smallskip
\hskip1cm
(b) $x_k=x_0-1$ and $x_{k-1}=1$ (an almost-culminating period).
\end{lem}
{\it Proof.} (i)  Recall that  if $p\ne 2$ is a prime,  then one and only one of the following
diophantine equations is solvable:
\begin{equation}\label{eq14}
\left\{
\begin{array}{ccc}
x^2-py^2 &=& -1,\\
x^2-py^2 &=& 2,\\
x^2-py^2 &=& -2,
\end{array}
\right.
\end{equation}
see \cite[p. 97 (Satz 3.21)]{P}.  Since  $p\equiv 3~\mod~4$,  one concludes
that $x^2-py^2=-1$ is not solvable \cite[p. 98 (Satz 3.23-24)]{P};  this happens if and only if $P=2k$ is
even (for otherwise the continued fraction of $\sqrt{p}$ would provide a
solution).

It is known,  that for even periods $P=2k$ the convergents  $A_i/B_i$ satisfy the
diophantine equation $A_{k-1}^2-pB_{k-1}^2=(-1)^k ~2$,  see \cite[p. 103]{P};
thus if $P\equiv 0~\mod~4$,  the equation $x^2-py^2=2$ is solvable and if
$P\equiv 2~\mod~4$,  then the equation $x^2-py^2=-2$ is solvable.
But  equation   $x^2-py^2=2$  (equation $x^2-py^2=-2$,  resp.) is solvable if and only if 
$p\equiv 7~\mod~8$  ($p\equiv 3~\mod~8$, resp.), see  \cite[Satz 3.23]{P}
(\cite[Satz 3.24]{P}, resp.)  Item (i) follows.

\medskip
(ii)   The equation $A_{k-1}^2-pB_{k-1}^2=(-1)^k ~2$ is a special case of 
equation $A_{k-1}^2-pB_{k-1}^2=(-1)^k ~Q_k$,  where $Q_k$ is the full
quotient of continued fraction \cite[p.92]{P};  therefore, $Q_k=2$.
One  can now apply  \cite[Satz 3.15]{P},  which says that for $P=2k$ and 
$Q_k=2$ the continued fraction of    $\sqrt{{\goth P}_{3~\mod~4}}$ is either
culminating (i.e.   $x_k=x_0$) or almost-culminating  (i.e. $x_k=x_0-1$ and $x_{k-1}=1$).
Lemma \ref{lm3} follows.
$\square$

\begin{lem}\label{lm4}
If $p\equiv 3~\mod~8$,   then  $c({\cal A}_{RM}^{(p,1)})=2$.
\end{lem}
{\it Proof.}  The proof proceeds by induction in period $P$, which is
in this case   $P\equiv 2~\mod~4$ by lemma \ref{lm3}.    We shall start
with $P=6$, since $P=2$ reduces to it,   see item (i) below.

\smallskip
(i)   Let $P=6$ be a culminating period;  then equation (\ref{eq5}) admits 
a general solution $[x_0,\overline{x_1,  2x_1,  x_0, ,  2x_1, x_1,  2x_0}]=\sqrt{x_0^2+4nx_1+2}$,
where $x_0=n(2x_1^2+1)+x_1$  \cite[p.101]{P}.   The solution depends on two
integer variables $x_1$ and $n$,  which is the maximal  possible number of  variables in this case;
therefore,  the dimension of the solution is $2$,  so as  complexity of the corresponding torus.   
Notice that the case  $P=2$ is obtained from $P=6$ by restriction to  $n=0$;  thus the complexity   
for $P=2$ is equal to $2$.  

\smallskip
(ii)   Let $P=6$ be an almost-culminating period;  
 then equation (\ref{eq5}) has 
a solution $[3s+1,\overline{2,  1,  3s ,  1,  2,  6s+2}]=\sqrt{(3s+1)^2+2s+1}$,
where $s$ is an integer variable  \cite[p. 103]{P}.  We encourage the
reader to verify,  that this solution is a restriction of solution (i) to $x_1=-1$
and $n=s+1$;   thus,  the dimension of our solution is $2$, so as the
complexity of the corresponding torus.     

\smallskip
(iii)  Suppose  a solution $[x_0,\overline{x_1,\dots, x_{k-1}, x_k, x_{k-1}, \dots, x_1, 2x_0}]$
with the (culminating or almost-culminating) period $P_0\equiv 3~\mod~8$ has dimension $2$;
let us show that a solution
\begin{equation}\label{eq15}
[x_0,\overline{y_1, x_1,\dots, x_{k-1}, y_{k-1},  x_k, y_{k-1},  x_{k-1}, \dots, x_1, y_1, 2x_0}]
\end{equation}
with period $P_0+4$ has also dimension $2$.  According to \cite{Web1}, 
if (\ref{eq15}) is a solution to the diophantine equation (\ref{eq5}),  then either
(i) $y_{k-1}=2y_1$ or (ii) $y_{k-1}=2y_1+1$ and $x_1=1$.  We proceed by showing 
that case (i)  is not possible for the square roots of prime numbers.   

Indeed, let to the contrary $y_{k-1}=2y_1$;  then the following system of equations
must be compatible:
\begin{equation}\label{eq16}
\left\{
\begin{array}{ccc}
A_{k-1}^2 &-&   pB_{k-1}^2=-2,\\
A_{k-1} &=& 2y_1A_{k-2}+A_{k-3},\\
B_{k-1} &=& 2y_1B_{k-2}+B_{k-3},
\end{array}
\right.
\end{equation}
 where $A_i, B_i$  are convergents and the first equation is solvable
 since $p\equiv 3~\mod~8$.  From the first equation, both convergents 
 $A_{k-1}$ and $B_{k-1}$ are odd numbers. (They are both odd or even,
 but the even must be excluded,  since $A_{k-1}$ and $B_{k-1}$ are relatively prime.)
 From the last two equations, the convergents $A_{k-3}$ and $B_{k-3}$
 are also odd.  Then the convergents $A_{k-2}$ and $B_{k-2}$ must
 be even,  since among six consequent  convergents $A_{k-1}, B_{k-1}, A_{k-2}, B_{k-2},
 A_{k-3}, B_{k-3}$ there are always two even;   but this is not possible,  because
 $A_{k-2}$ and $B_{k-2}$ are relatively prime. Thus, $y_{k-1}\ne 2y_1$.

 Therefore (\ref{eq15}) is a solution of the diophantine equation (\ref{eq5})
 if and only if  $y_{k-1}=2y_1+1$ and $x_1=1$;  the dimension of such a solution
 coincides with the dimension of solution 
 $[x_0,\overline{x_1,\dots, x_{k-1}, x_k, x_{k-1}, \dots, x_1, 2x_0}]$,
 since for two new integer variables $y_1$ and $y_{k-1}$ one gets
 two new constraints.  Thus,  the dimension of solution (\ref{eq15})
 is $2$,  so as the  complexity of the corresponding torus. 
 Lemma \ref{lm4} follows.
 $\square$

\begin{lem}\label{lm5}
If $p\equiv 7~\mod~8$,   then  $c({\cal A}_{RM}^{(p,1)})=1$.
\end{lem}
{\it Proof.} 
 The proof proceeds by induction in period $P\equiv 0~\mod~4$, see 
  lemma \ref{lm3};  we  start with $P=4$. 
  
\smallskip
(i)   Let $P=4$ be a culminating period;  then equation (\ref{eq5}) admits 
a solution $[x_0,\overline{x_1,x_2,x_1,2x_0}]= \sqrt{x_0^2+m(x_1x_2+1)-x_2^2}$,
where $x_2=x_0$,   see  example  \ref{ex1}  for the details.   
Since  the polynomial  $m(x_0x_1+1)$  under   the square root
represents   a prime number,  we have   $m=1$;   the latter equation
is not solvable in integers $x_0$ and $x_1$,   since 
$m=x_0(x_0x_1+3)x_1^{-1}(x_0x_1+2)^{-1}$.   Thus,  there are no solutions
of (\ref{eq5})  with the culminating period $P=4$.

\smallskip
(ii)   Let $P=4$ be an almost-culminating period;  then equation (\ref{eq5}) admits 
a solution  $[x_0,\overline{1, x_0-1, 1, 2x_0}]= \sqrt{(x_0+1)^2-2}$.
The dimension of this solution was proved to be $1$,  see example \ref{ex1};
thus,  the complexity of the corresponding torus is equal  to $1$.

\smallskip
(iii)  Suppose  a solution $[x_0,\overline{x_1,\dots, x_{k-1}, x_k, x_{k-1}, \dots, x_1, 2x_0}]$
with the (culminating or almost-culminating) period $P_0\equiv 7~\mod~8$ has dimension $d=1$.
It can be shown by the same argument as in lemma \ref{lm4},  that for 
 a solution of the form (\ref{eq15}) having the  period $P_0+4$ the  dimension remains
 the same, i.e. is equal to $1$;  we leave details to the reader. Thus, complexity of the corresponding
 torus is equal to $1$. Lemma \ref{lm5} follows.
 $\square$

\begin{lem}\label{lm6}
{\bf (\cite[p. 78]{G})}
\begin{equation}\label{eq17}
 rk_{\Bbb Q} ~({\cal E}_{CM}^{(-p,1)})=
 \cases{1, & if  ~$p\equiv 3~\mod~8$\cr
              0, & if   ~$p\equiv 7~\mod~8.$}
\end{equation}
\end{lem}
\begin{rmk}
\textnormal{
For the sake of clarity, notice that \cite[Theorem 22.4.2]{G} for $p\equiv 3~\mod~8$ 
gives an estimate $rk_{\Bbb Q} ~({\cal E}_{CM}^{(-p,1)})\le 1$ only. 
However, on the next page \cite[p. 79]{G} it is remarked that the actual 
value $rk_{\Bbb Q} ~({\cal E}_{CM}^{(-p,1)})=1$. 
The same is true of the condition $Cl(L)_p\left({3p-1\over 4}\right)=(0)$,
which can be omitted for  the value $rk_{\Bbb Q} ~({\cal E}_{CM}^{(-p,1)})=1$,
{\it ibid.}
}  
\end{rmk}

\bigskip
Comparing lemmas \ref{lm4} and \ref{lm5} with the formula (\ref{eq17}),  one gets the conclusion of theorem \ref{thm1}.
$\square$

\section{Examples}
To illustrate theorem \ref{thm1},  let us consider the  ${\Bbb Q}$-curves ${\cal E}_{CM}^{(-p,1)}$,
where  $p<100$.  Our results are shown in Figure 1.

\begin{figure}
\begin{tabular}{c|c|c|c}
\hline
&&&\\
$p\equiv 3~\mod~4$ & $rk_{\Bbb Q}({\cal E}_{CM}^{(-p,1)})$ & $\sqrt{p}$ & $c({\cal A}_{RM}^{(p,1)})$\\
&&&\\
\hline
$3$ & $1$ & $[1,\overline{1,2}]$ & $2$\\
\hline
$7$ & $0$ & $[2,\overline{1,1,1,4}]$ & $1$\\
\hline
$11$ & $1$ & $[3,\overline{3,6}]$ & $2$\\
\hline
$19$ & $1$ & $[4,\overline{2,1,3,1,2,8}]$ & $2$\\
\hline
$23$ & $0$ & $[4,\overline{1,3,1,8}]$ & $1$\\
\hline
$31$ & $0$ & $[5,\overline{1,1,3,5,3,1,1,10}]$ & $1$\\
\hline
$43$ & $1$ & $[6,\overline{1,1,3,1,5,1,3,1,1,12}]$ & $2$\\
\hline
$47$ & $0$ & $[6,\overline{1,5,1,12}]$ & $1$\\
\hline
$59$ & $1$ & $[7,\overline{1,2,7,2,1,14}]$ & $2$\\
\hline
$67$ & $1$ & $[8,\overline{5,2,1,1,7,1,1,2,5,16}]$ & $2$\\
\hline
$71$ & $0$ & $[8,\overline{2,2,1,7,1,2,2,16}]$ & $1$\\
\hline
$79$ & $0$ & $[8,\overline{1,7,1,16}]$ & $1$\\
\hline
$83$ & $1$ & $[9,\overline{9,18}]$ & $2$\\
\hline
\end{tabular}

\caption{The ${\Bbb Q}$-curves ${\cal E}_{CM}^{(-p,1)}$  with  $p<100$.}
\end{figure}




\vskip1cm
\textsc{Department of Mathematics and Computer Science, St.~John's University, 8000 Utopia Parkway,  
New York,  NY 11439, United States;} ~\textsc{E-mail:} {\sf igor.v.nikolaev@gmail.com}


\begin{thebibliography}{100}
\bibitem{G}
B.~H.~Gross, Arithmetic on Elliptic Curves with Complex Multiplication,
Lecture Notes Math. 776 (1980), Springer. 



\bibitem{Man1}
Yu.~I.~Manin, Real multiplication and noncommutative geometry,
in ``Legacy of Niels Hendrik Abel'', 685-727, Springer, 2004. 


\bibitem{M}
G.~J.~Murphy,  $C^*$-Algebras and Operator Theory, Academic Press, 1990.





\bibitem{N}
I.~Nikolaev, Noncommutative Geometry,
De Gruyter Studies in Math. 66,   Berlin,  2017.


\bibitem{Nik3}
I.~Nikolaev, Noncommutative geometry of twists, 
arXiv:1712.07516

\bibitem{P}
O.~Perron, Die Lehre von den Kettenbr\"uchen,  Bd.1,  Teubner, 1954.  


\bibitem{Rie2}
M.~A.~Rieffel, Non-commutative tori -- a case study of non-commutative
differentiable manifolds,  Contemp. Math. 105 (1990), 191-211. 


\bibitem{S}
J.~H.~Silverman, Advanced Topics in the Arithmetic of Elliptic Curves,
GTM 151, Springer 1994.



\bibitem{StaVdb1}
J.~T.~Stafford and M.~van ~den ~Bergh, Noncommutative curves and noncommutative
surfaces, Bull. Amer. Math. Soc. 38 (2001), 171-216. 



\bibitem{Web1}
K.~Weber,  Kettenbr\"uche  mit  kulminierenden  und  fastkuminierenden   
Perioden,  Sitzungsber. der Bayer. Akademie d. Wissenschaften zu M\"unchen,
mathemat.-naturwissen. Abteilung  (1926), 41-62.  


\end{thebibliography}
\end{document}